\documentclass{amsart}
\usepackage[utf8]{inputenc}
\usepackage[whole]{bxcjkjatype}

\usepackage{bm,amsfonts}
\usepackage{amsmath}
\usepackage{amssymb}
\usepackage{graphicx}
\usepackage{hyperref}
\usepackage{physics}
\usepackage{color}
\usepackage{tikz}

\newtheorem*{theorem}{Theorem}
\newtheorem*{conjecture}{Conjecture}

\date{\today}

\title{A refinement of Sato--Tate conjecture}
\author{Taro Kimura \\ 木村太郎}
\date{\today}
\address{Institut de Math\'ematiques de Bourgogne, Universit\'e Bourgogne France-Comt\'e}

\begin{document}

\maketitle

\begin{abstract}
    We propose a refined version of the Sato--Tate conjecture about the spacing distribution of the angle determined for each prime number.
    We also discuss its implications on $L$-function associated with elliptic curves in the relation to random matrix theory.
\end{abstract}


\section{Introduction}

The Sato--Tate (ST) conjecture, independently proposed by M.~Sato and J.~Tate~\cite{Tate:1965}, then proved by Clozel--Harris--Shepherd-Barron--Taylor~\cite{Clozel:2008PMIHES,Taylor:2008PMIHES,Harris:2010AM} after 40 years, is about a statistical property of the number of points on the elliptic curve over the finite field.
Let $\#E(\mathbb{F}_p)$ be the number of points on the elliptic curve $E$ over the finite field $\mathbb{F}_p$ with $p$ a prime number, and $N$ the conductor of $E$.
For each pair $(p,E)$, we define $a_p = p + 1 - \#E(\mathbb{F}_p)$ for $p \nmid N$ and $a_p = p  - \#E(\mathbb{F}_p)$ for $p \mid N$.
Hasse's theorem (also equivalent to the Ramanujan conjecture in some cases) on elliptic curves states that
\begin{align}
    \abs{\frac{a_p}{2 \sqrt{p}}} \le 1
    \, ,
    \label{eq:Hasse_thm}
\end{align}
from which we define the angle $\theta_p \in [0,\pi]$ for a pair $(E,p)$,
\begin{align}
    \cos \theta_p = \frac{a_p}{2 \sqrt{p}}
    \, .
    \label{eq:theta_def}
\end{align}

\begin{theorem}[The Sato--Tate conjecture. Proof by Clozel--Harris--Shepherd-Barron--Taylor]
Suppose that the elliptic curve $E$ has no complex multiplication.
Then, the probability such that the angle $\theta_p$ is found in the interval $[a,b] \subseteq [0,\pi]$ is given by
\begin{align}
    \mathbb{P}\qty[ a \le \theta_p \le b ] 
    := \lim_{X \to \infty} \mathbb{P}_X\qty[ a \le \theta_p \le b ] 
    = \int_a^b \dd{\theta} \rho(\theta)
    \, ,
    \label{eq:1-pt}
\end{align} 
where
\begin{align}
    \mathbb{P}_X\qty[ a \le \theta_p \le b ] 
    = \frac{\#\qty{ a \le \theta_p \le b , p \le X}}{\#\qty{ p \le X}}
\end{align}
with the density function given by
\begin{align}
    \rho(\theta) = \frac{2}{\pi} \sin^2 \theta
    \, .
    \label{eq:density_fn}
\end{align}
\end{theorem}

In fact, the density function discussed here is interpreted as a one-point correlation function of the angle distributed on $[0,\pi]$.
We may obtain it from the probability~\eqref{eq:1-pt} in the following limit,
\begin{align}
    \lim_{\dd{\theta} \to \varepsilon \ll 1}
    \mathbb{P} \qty[\theta \le \theta_p \le \theta + \dd{\theta}]
    = \rho(\theta) \dd{\theta}
    \, .
\end{align}
From this point of view, it is natural to consider more generic statistical quantities to characterize the distribution of the angles.
For example, we may consider the $n$-point generalization of the probability~\eqref{eq:1-pt} for $0 \le a_i \le b_i \le \pi$, $i = 1,\ldots,n$, \begin{align}
    \mathbb{P} \qty[ a_i \le \theta_{p_i} \le b_i, i = 1, \ldots,n]
    & := \lim_{X \to \infty} \mathbb{P}_X \qty[ a_i \le \theta_{p_i} \le b_i, i = 1, \ldots,n]
    \nonumber \\
    & = 
    \lim_{X \to \infty} \frac{\#\qty{ a_i \le \theta_{p_i} \le b_i , i = 1, \ldots, n, p_i \le X}}{\#\qty{ p \le X}}
    \, .
\end{align}
Then, the $n$-point function is similarly obtained through the limit
\begin{align}
    \lim_{\dd{\theta}_i \to \varepsilon \ll 1} 
    \mathbb{P} \qty[ \theta_i \le \theta_{p_i} \le \theta_i + \dd{\theta}_i, i = 1, \ldots,n]
    =
    \rho_n(\theta_1,\ldots,\theta_n) \dd{\theta}_1 \cdots \dd{\theta}_n
    \, .
    \label{eq:n-pt}
\end{align}

In this paper, we do not directly consider such $n$-point functions, but instead study another statistical quantity about spacings between the angles.
We define an integral map from $[0,\pi]$ to $[0,1]$,
\begin{align}
    \Theta (\theta) = \int_0^{\theta} \dd{\theta'} \rho(\theta')
    \, .
\end{align}
Since $\dd{\Theta} = \rho(\theta) \dd{\theta}$, the variable $\Theta$ is uniformly distributed on $[0,1]$, while the variable $\theta$ is distributed on $[0,\pi]$ with the density profile $\rho(\theta)$.
We introduce a set of variables $(\Theta(\theta_p))_{p \in \text{primes}} = (\Theta_1 \le \Theta_2 \le \cdots )$, that we call unfolded variables.
Then, we have the following conjecture based on numerical computations of the unfolded variables shown in Section~\ref{sec:numerics}.
\begin{conjecture}[Spacing distribution of the unfolded angles]
\label{conj:sp}
For the elliptic curve $E$, for which the ST conjecture holds, the next$^k$\! nearest neighbor spacing distribution%
\footnote{%
Nearest neighbor $(k=0)$, Next nearest neighbor $(k=1)$, Next-next nearest neighbor $(k=2)$, etc.
}%
of the unfolded variables is given by the Poisson distribution,
\begin{align}
    p_{k}(s) := 
    \lim_{X \to \infty}
    \mathbb{P}_X \qty[ \Theta_{i + k + 1} - \Theta_i = \frac{s}{\#\{p \le X\} } ]
    = \frac{s^k \, e^{-s}}{k!}
    \, .
\end{align}
The spacing distribution is normalized
\begin{align}
    \int_0^\infty \dd{s} p_{k}(s) = 1
    \, , \qquad
    \int_0^\infty \dd{s} s \, p_{k}(s) = k+1
    \, .
    \label{eq:sp_normalization}
\end{align}
\end{conjecture}
This is a refined version of the ST conjecture in a sense that it concerns more detailed statistical property of the angles beyond the density function addressed in the ordinary ST conjecture.

The Poisson distribution appears as a spacing distribution of uncorrelated random variables.
See, for example, \cite{Mehta:2004RMT}.
Therefore, this conjecture implies that the angles $\theta_p$ and the unfolded variables $\Theta_i$ statistically behave as uncorrelated random variables.
From this point of view, it is speculated that the $n$-point function~\eqref{eq:n-pt} would be simply given as a product of the density functions,
\begin{align}
        \rho_n(\theta_1,\ldots,\theta_n) = \prod_{i=1}^n \rho(\theta_i)
        \, .
    \end{align}
We remark that a similar behavior is observed for prime numbers~\cite{Liboff:1998IJTP,Timberlake:2006AJP,Timberlake:2007,Wolf:2014PRE}:
The spacing distribution of prime numbers is given by the Poisson distribution.
Compared to these results, one can more clearly see that the spacing distribution of the angles is given by the Poisson distribution since the angles are distributed on the finite interval $[0,\pi]$, whereas the prime numbers are distributed on the infinite interval $[0,\infty]$.

\subsubsection*{Acknowledgments}

I would like to thank S.~Koyama and N.~Kurokawa for valuable discussions and useful comments on the draft.
In fact, this study was motivated by Kurokawa-sensei's seminar explaining a similarity between the Sato--Tate conjecture and Wigner's semi-circle law of random matrices.
I'm grateful for his exposition.
This work has been supported in part by ``Investissements d'Avenir'' program, Project ISITE-BFC (No.~ANR-15-IDEX-0003), and EIPHI Graduate School (No.~ANR-17-EURE-0002).

\section{$L$-function for elliptic curves}

We address an implication of the refined ST conjecture in the context of the $L$-function associated with elliptic curves.

\subsection{Global zero and local zero}

In addition to the \emph{algebraic} formulation shown above, the numbers $a_p$ also appear on the \emph{analytic} side as the Fourier coefficients of the cusp form associated with the elliptic curve $E$ (See, for example, \cite{Koblitz:1993}),
\begin{align}
    f(q) = \sum_{n = 1}^\infty a_n \, q^n
    \, .
\end{align}
We define the $L$-function associated with the elliptic curve $E$, a.k.a., the Hasse--Weil $L$-function, from the Fourier coefficients of the cusp form, which has the Euler product form due to the multiplicative property of $a_p$'s,%
\begin{align}
    L(s;E) 
    & = \sum_{n = 1}^\infty \frac{a_n}{n^s} 
    \nonumber \\ &
    = \prod_{p \nmid N} \qty(1 - a_p \, p^{-s} + p^{1 - 2s})^{-1} \prod_{p \mid N} \qty(1 - p^{-s})^{-1}
    \, .
    \label{eq:L_fn}
\end{align}
We then define the local zeta function for a pair $(p,E)$, which shows the following factorization,
\begin{align}
    Z_p(s;E) 
    = 1 - a_p \, p^{-s} + p^{1 - 2s} 
    = (1 - e^{+ i \theta} p^{1/2 - s}) (1 - e^{- i \theta} p^{1/2 - s})
    \, .
    \label{eq:Z_fn}
\end{align}
In other words, the relation~\eqref{eq:Hasse_thm} is equivalent to the local version of Riemann Hypothesis:
The zeros of the local zeta function, that we call the local zeros, are found on the critical line $\Re(s) = \frac{1}{2}$.
This is for a cusp form of weight two, and a similar formulation is also possible for higher weight cases.%
\footnote{%
The local zeta function for the cusp form of weight $k$ is $\displaystyle Z_p(s;E) = 1 - a_p \, p^{-s} + p^{k - 1 - 2s}$, so that the critical line is given by $\Re(s) = (k-1)/2$.
In Section~\ref{sec:numerics}, we will discuss the modular discriminant, which is a cusp form of weight 12.
}

From this point of view, the ST conjecture about the distribution of the angles describes how the imaginary part of the local zero is distributed on the critical line.
In particular, our conjecture implies that the imaginary part of the local zero behaves as an uncorrelated random number.
This is in contrast to the statistical behavior of zeros of the $L$-function, that we instead call the global zeros.
It has been known that the statistical property of the nontrivial zeros of the zeta function agrees with eigenvalue statistics of random matrices~\cite{Montgomery:1973,Odlyzko:1987MC}.
This is true also for the $L$-function associated with elliptic curves~\cite{Bourgade:2013}.
Then, we have a question: 
What is a counterpart of the local zero in random matrix theory?
In the context of Gaussian random matrix theory, we typically consider statistical properties of eigenvalues exhibiting nontrivial correlation, rather than uncorrelated matrix elements which are independently distributed.
Namely, eigenvalues are correlated, whereas matrix elements are not.
This observation leads to the following correspondence:
    \begin{center}
        \begin{tabular*}{.7\textwidth}{@{\extracolsep{\fill}}ccc}\hline\hline
            Correlation & Random matrices & Zeta function \\ \hline
            Correlated & Eigenvalue & Global zero \\
            Uncorrelated & Matrix element & Local zero \\
            \hline\hline
        \end{tabular*}
    \end{center}
Recalling the Euler product of the $L$-function~\eqref{eq:L_fn}, the local zeta function is a building block of the total $L$-function.
A similar interpretation is possible for the ordinary zeta function in the Euler product form for $\Re(s) > 1$,
\begin{align}
    \zeta(s) = \sum_{n=1}^\infty \frac{1}{n^s} = \prod_{p \in \text{primes}} \qty(1 - p^{-s})^{-1}
    \, ,
\end{align}
because the spacing distribution of prime numbers is also the Poisson distribution as mentioned earlier.
This seems analogous to the relation between matrix elements and eigenvalues of random matrices:
Correlated variables are constructed from uncorrelated variables.

\subsection{Group structure and higher genus/rank generalization}

We remark that, introducing $U(\theta) = \operatorname{diag}\qty(e^{+i\theta},e^{-i\theta}) \in \mathrm{SU}(2)$, we may write the local zeta function~\eqref{eq:Z_fn} as a characteristic polynomial,
\begin{align}
    Z_p(s;E) = \det \qty(1 - U(\theta_p) p^{1/2-s} )
    \, .
\end{align}
The local zeta function depends only on the conjugacy class parametrized by the angle $\theta_p$ since it is invariant under the SU(2) transform, $U(\theta_p) \to g^{-1} U(\theta_p) g$ for ${}^\forall g \in \mathrm{SU}(2)$.
In fact, the measure $\rho(\theta) \dd{\theta}$ agrees with the Haar measure on SU(2).
This factor is also interpreted as follows:
Introducing another variable, $x = \cos x$, the measure is given by $\rho(\theta) \dd{\theta} = \frac{2}{\pi} \sqrt{1 - x^2} \dd{x}$, which is obtained by the projection onto the interval $x \in [-1,1]$ from a three sphere $\mathbb{S}^3 = \{ x^2 + y^2 + z^2 + w^2 = 1 \}$, identified with SU(2).
From this point of view, our conjecture implies that we have an uncorrelated homogeneous map from prime numbers to SU(2).

A similar interpretation may be possible for higher genus cases (hyperelliptic curves).
It has been proposed that a higher rank group $G = \mathrm{USp}(2g)$ corresponds to genus $g$ curves, which is reduced to SU(2) for $g = 1$ elliptic curves~\cite{Katz:1999}.
A naive expectation is that we would similarly have uncorrelated homogeneous map from prime numbers to $G$.
However, in this case, we should use $g$ angles for each pair $(p,E)$, $(\theta_p^{(i)})_{i=1,\ldots,g}$, so that statistical behaviors of the angles would be more involved because of a correlation among them due to the Haar measure on $G$.
Further examination is required in this case.

\section{Numerical computations}\label{sec:numerics}

We demonstrate the statistical property of the angles with numerical computations.
We consider the following cusp forms as examples:
\begin{center}
    \begin{tabular*}{.8\textwidth}{@{\extracolsep{\fill}}rcccc}\hline\hline
        & \multicolumn{2}{c}{$f(q)$} & $N$ & $\#\{ p \le X \}$ \\ \hline
        (a) & $\displaystyle \eta(q)^2 \eta(q^{11})^2$ & ($E$: $y^2 + y = x^3 - x^2$) & 11 & 2000 \\
        (b) & $\displaystyle \eta(q^2)^2 \eta(q^{10})^2$ & ($E$: $y^2 = x^3 + x^2 - x$) & 2, 5 & 2000 \\
        (c) & $\displaystyle \eta(q)^{24}$ & (modular discriminant $\Delta$) & $-$ & 10000
        \\
        \hline\hline
    \end{tabular*}
\end{center}
where the Dedekind eta function is denoted by
\begin{align}
    \eta(q) = q^{1/24} \prod_{n=1}^\infty (1 - q^n)
    \, .
\end{align}
The right most column shows the number of primes that we use in the computations.
We remark that, since the weight of the modular discriminant (example (c)) is 12, the angle is defined as $\cos \theta_p = a_p / 2 p^{11/2}$, instead of the previous definition \eqref{eq:theta_def}.
In this case, the coefficient $a_p$ is also denoted by $\tau(p)$, which is called the Ramanujan tau function.
Similarly, the local zeta function is given by $\displaystyle Z_p(s,\Delta) = 1 - a_p \, p^{-s} + p^{11-2s}$.

\begin{figure}
    \begin{center}
    \begin{tikzpicture}
    
    \node at (0,0) {\includegraphics[height=3.5cm]{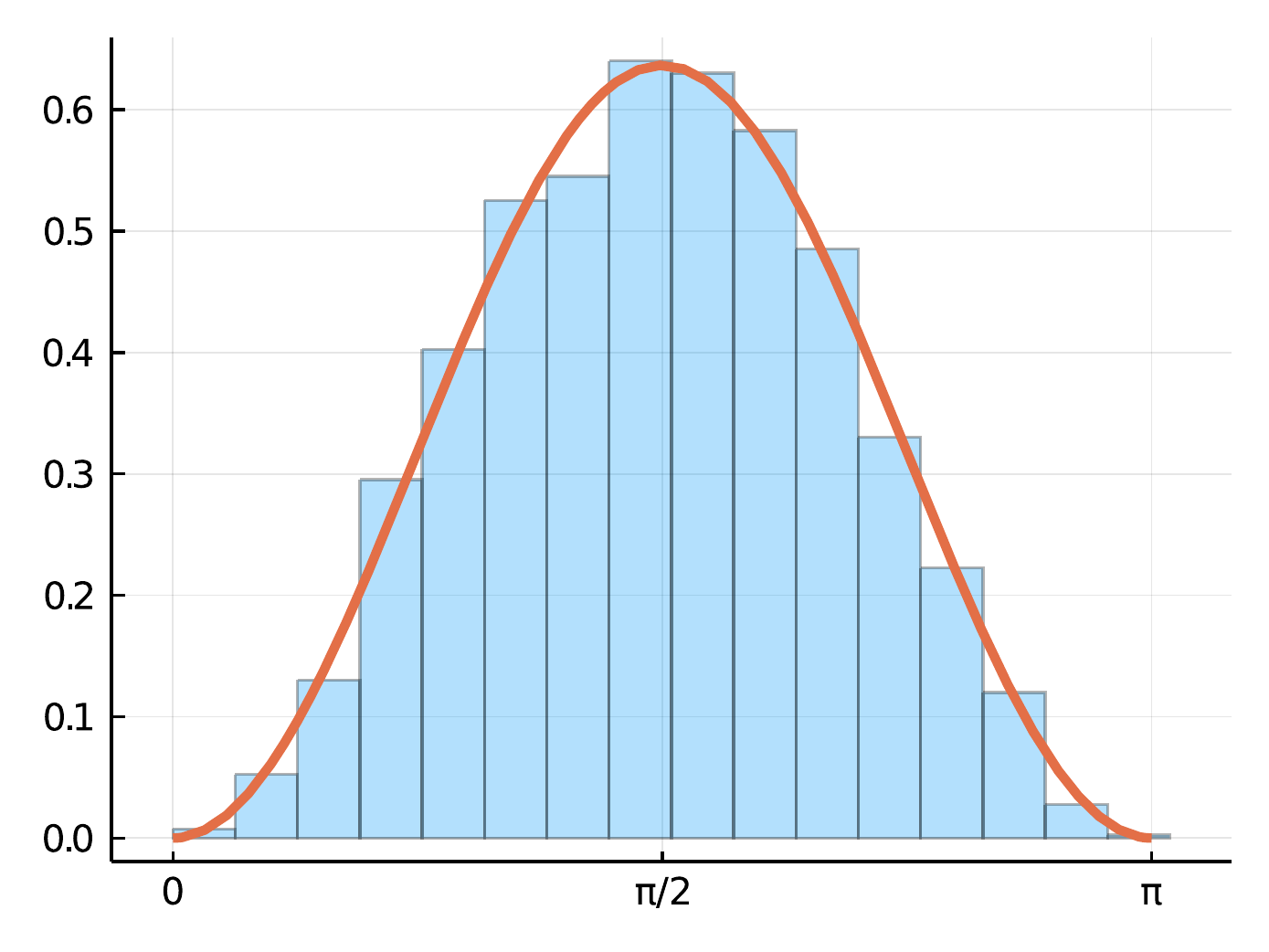}};
    \node at (1.8,1.8) {(a)};
    \node at (2.2,-1.2) {$\theta$};
    \node at (-1.9,1.9) {$\rho(\theta)$};
    
    \begin{scope}[shift = {(5.5,0)}]
    \node at (0,0) {\includegraphics[height=3.5cm]{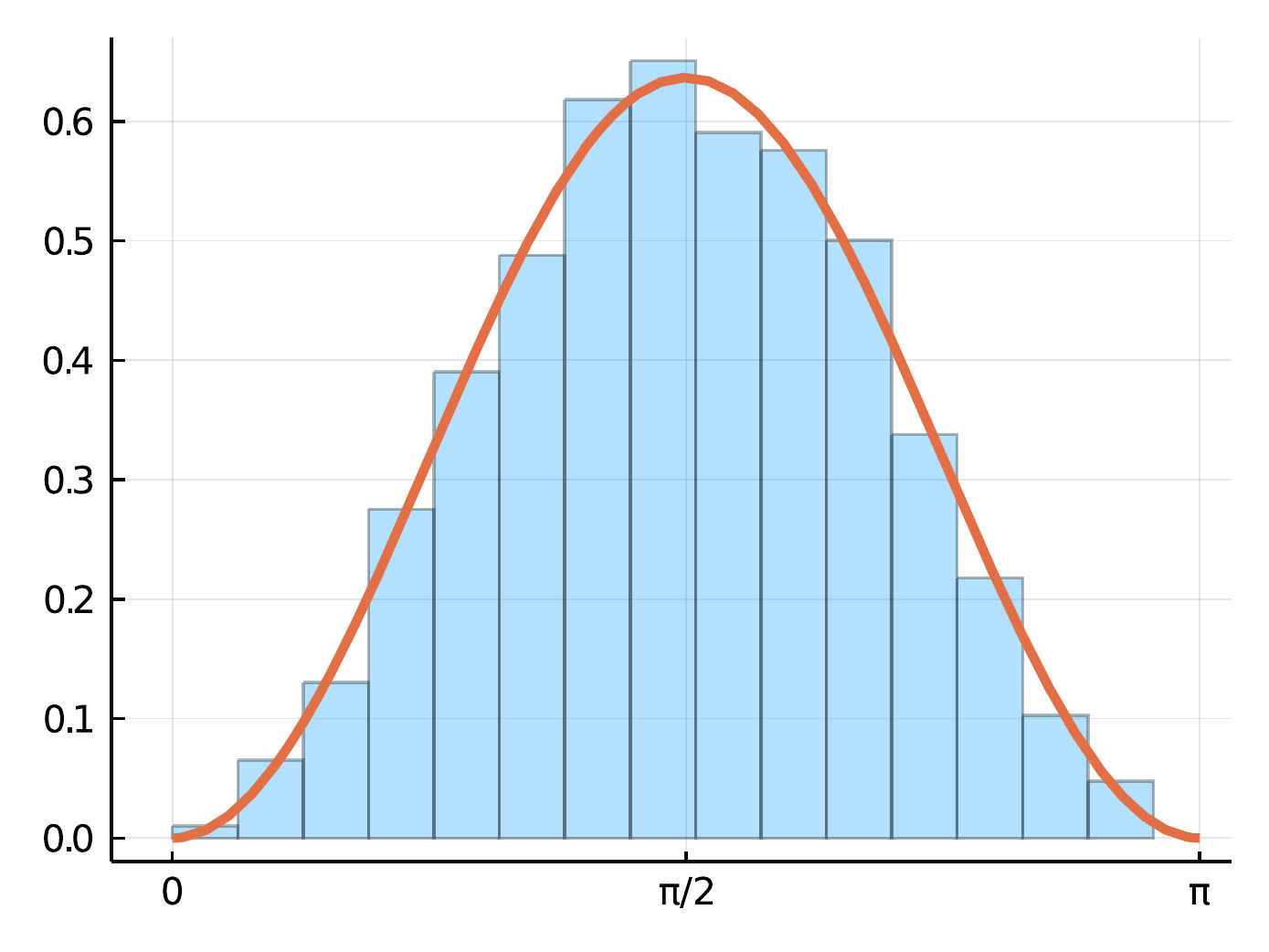}};
    \node at (1.8,1.8) {(b)};
    \node at (2.2,-1.2) {$\theta$};
    \node at (-1.9,1.9) {$\rho(\theta)$};
    \end{scope}
    
    \begin{scope}[shift = {(11,0)}]
    \node at (0,0) {\includegraphics[height=3.5cm]{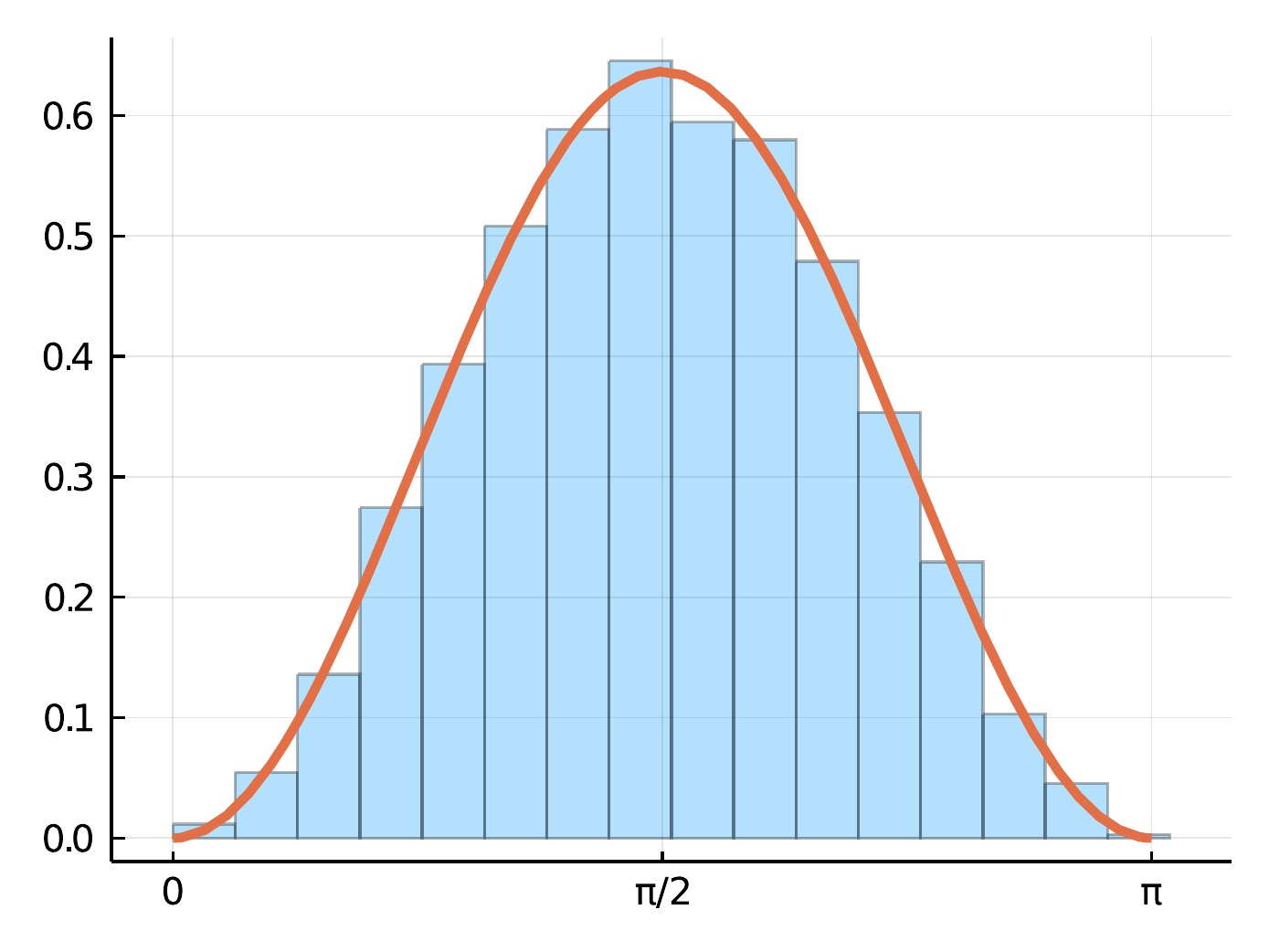}};
    \node at (1.8,1.8) {(c)};
    \node at (2.2,-1.2) {$\theta$};
    \node at (-1.9,1.9) {$\rho(\theta)$};
    \end{scope}
    
    \end{tikzpicture}
    \end{center}
    \caption{Numerical computation of the density function $\rho(\theta)$. 
    The red line shows the curve $\displaystyle \frac{2}{\pi} \sin^2 \theta$.}
    \label{fig:rho}
\end{figure}
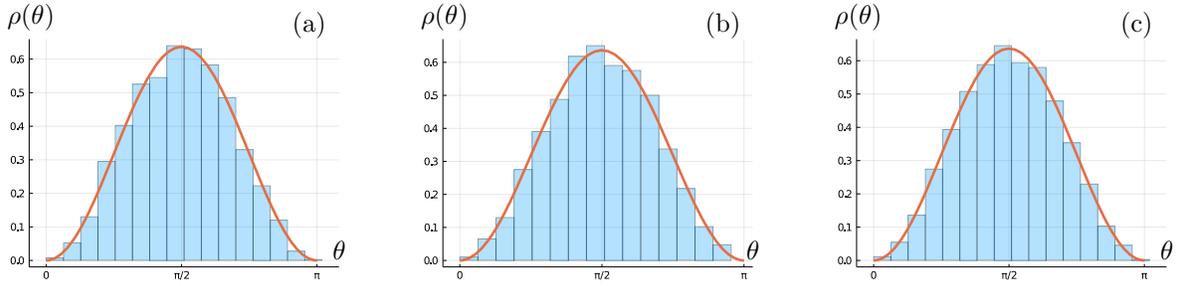

\begin{figure}
    \begin{center}
    \begin{tikzpicture}
    \node at (0,0) {\includegraphics[height=3.5cm]{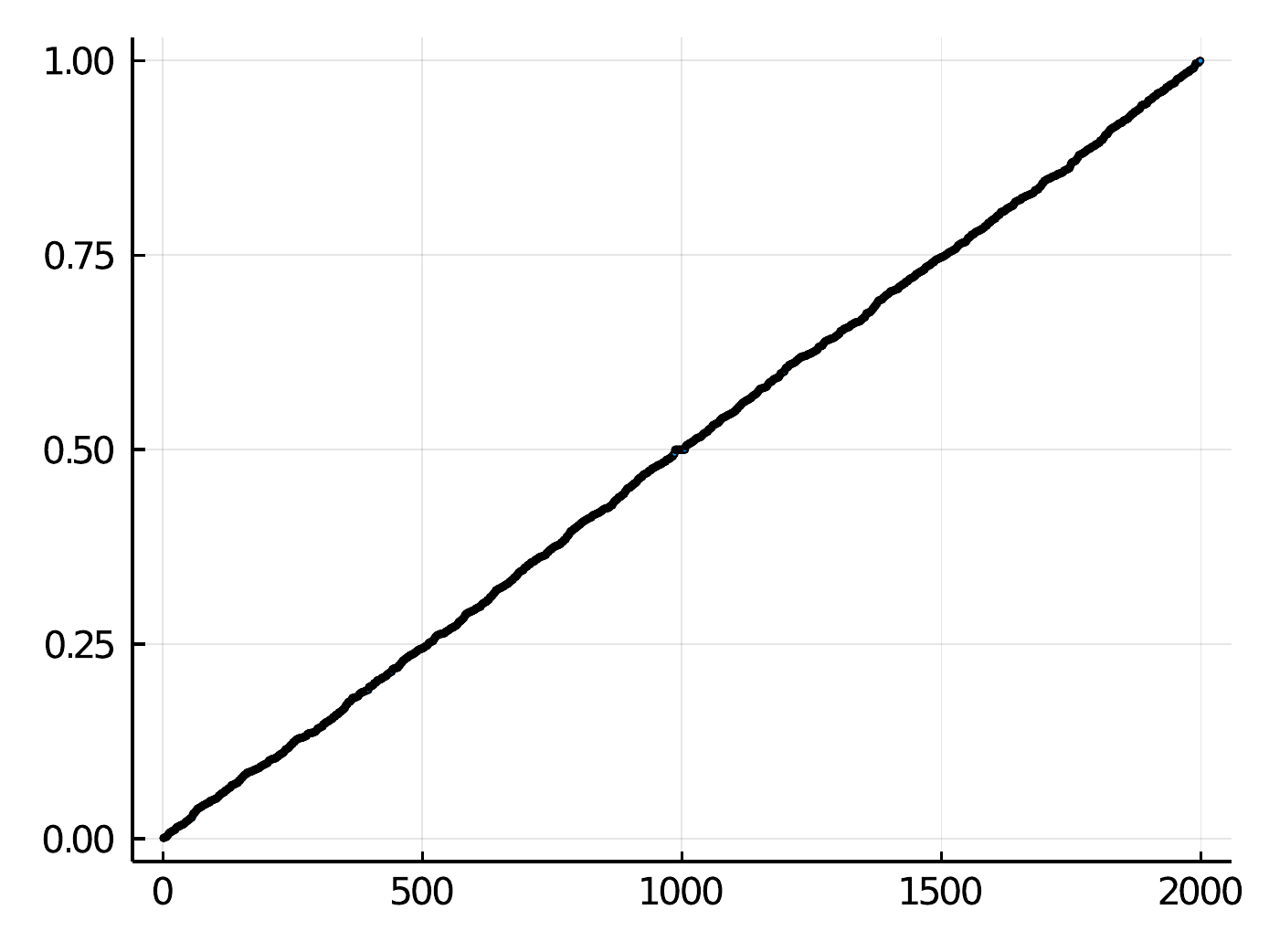}};
    \node at (1.8,1.8) {(a)};
    \node at (-1.9,1.9) {$\Theta_i$};
    \node at (2.3,-1.2) {$i$};
    
    \begin{scope}[shift = {(5.5,0)}]
    \node at (0,0) {\includegraphics[height=3.5cm]{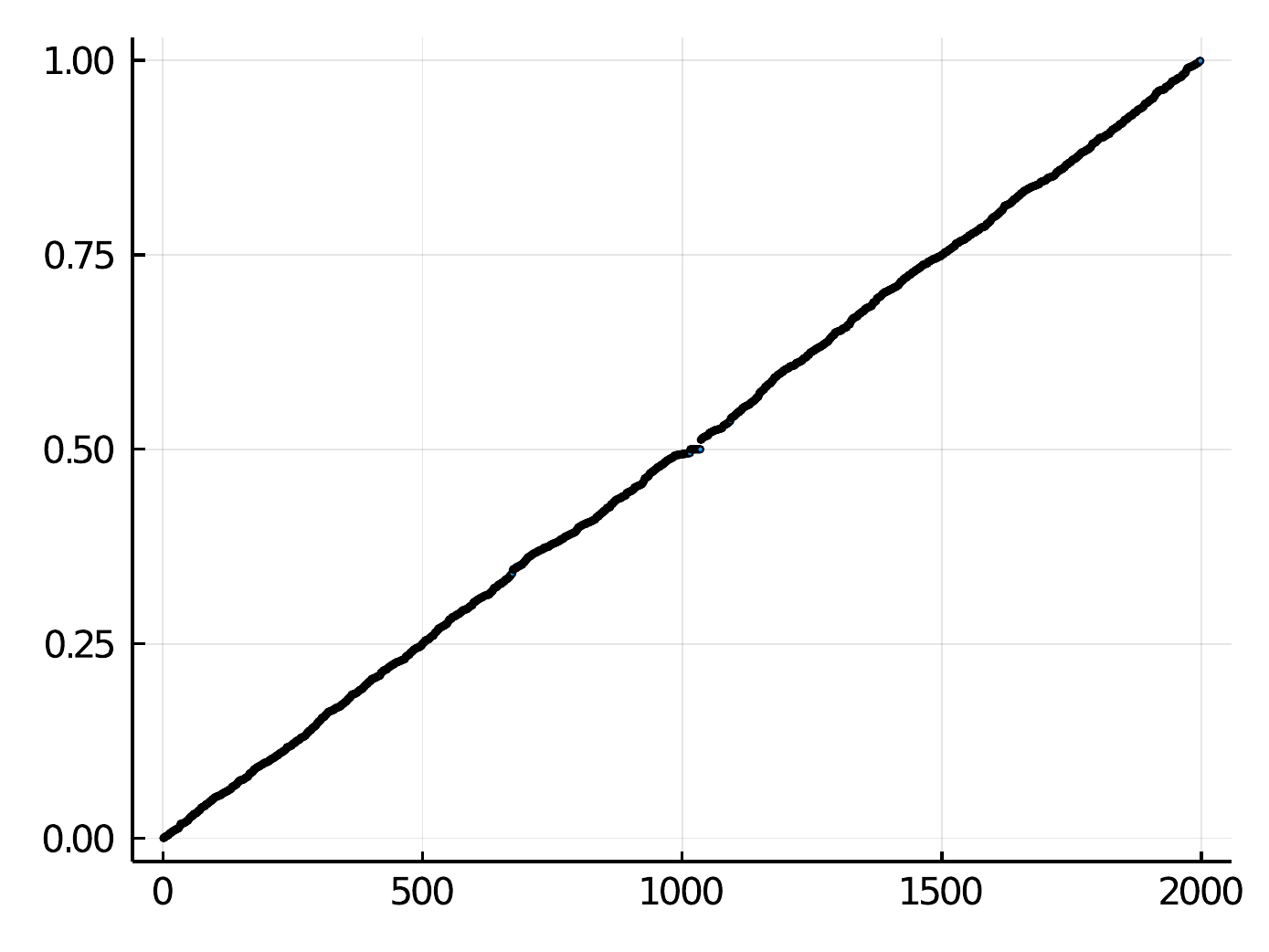}};
    \node at (1.8,1.8) {(b)};
    \node at (-1.9,1.9) {$\Theta_i$};
    \node at (2.3,-1.2) {$i$};
    \end{scope}
    
    \begin{scope}[shift = {(11,0)}]
    \node at (0,0) {\includegraphics[height=3.5cm]{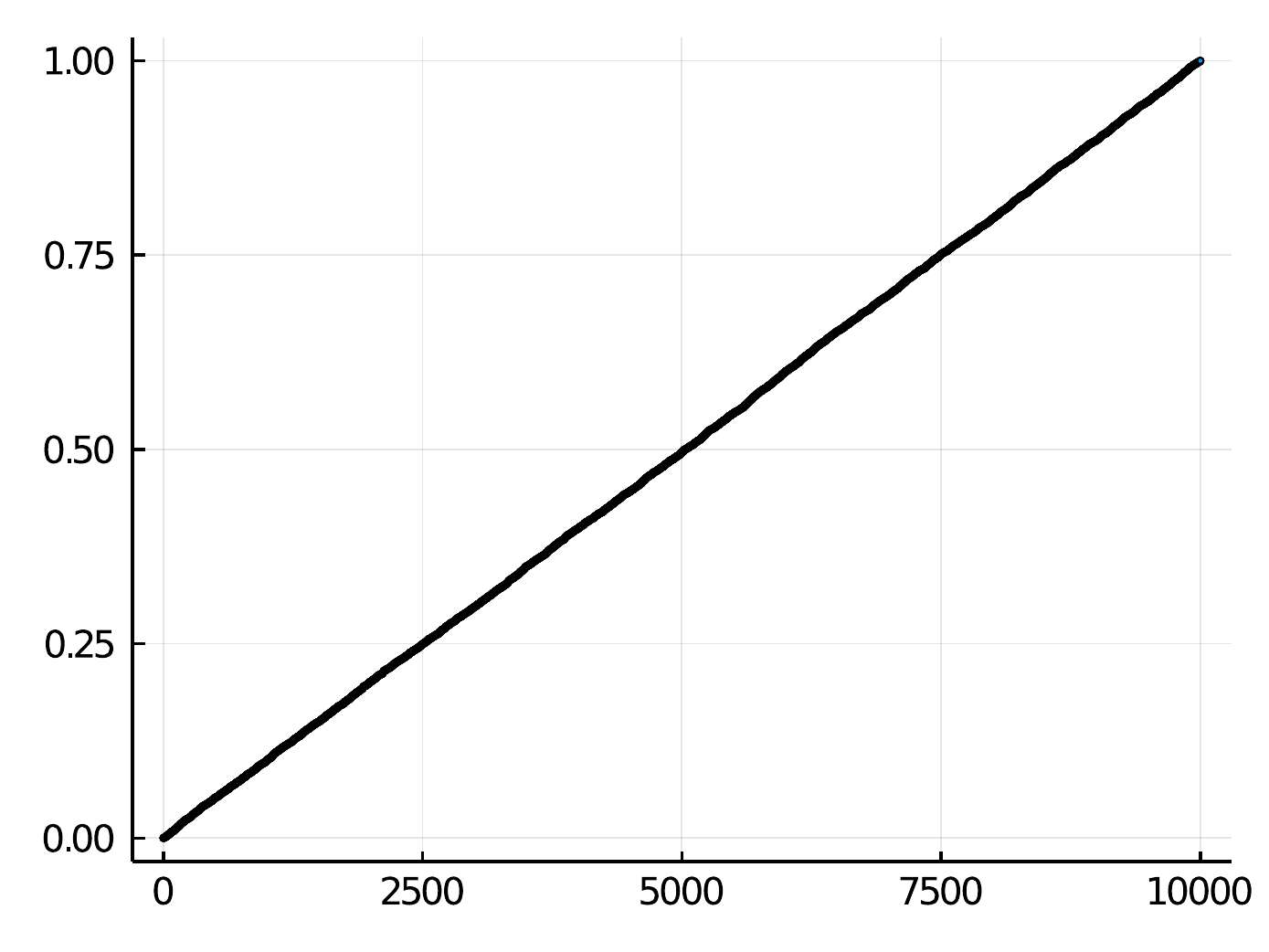}};
    \node at (1.8,1.8) {(c)};
    \node at (-1.9,1.9) {$\Theta_i$};
    \node at (2.3,-1.2) {$i$};
    \end{scope}
    \end{tikzpicture}
    \end{center}
    \caption{The unfolded variables $\Theta_i$ for $i = 1, \ldots, \#\{p \le X\}$.}
    \label{fig:unfold}
\end{figure}
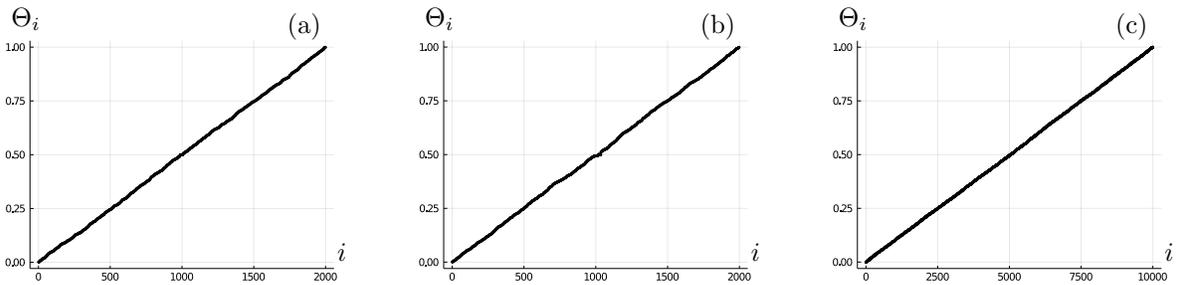

Fig.~\ref{fig:rho} shows numerical computation of the density function of the angles $\theta_p$ for three examples.
We see that it agrees with the curve $\displaystyle \rho(\theta) = \frac{2}{\pi} \sin^2 \theta$ as stated by the ST conjecture.
Fig.~\ref{fig:unfold} shows the unfolded variables $(\Theta_i)_{i = 1,\ldots,\#\{p \le X\}}$, which are uniformly distributed on the interval $[0,1]$.

\begin{figure}
    \begin{center}
        
    \begin{tikzpicture}
    \node at (0,0) {\includegraphics[height=3.5cm]{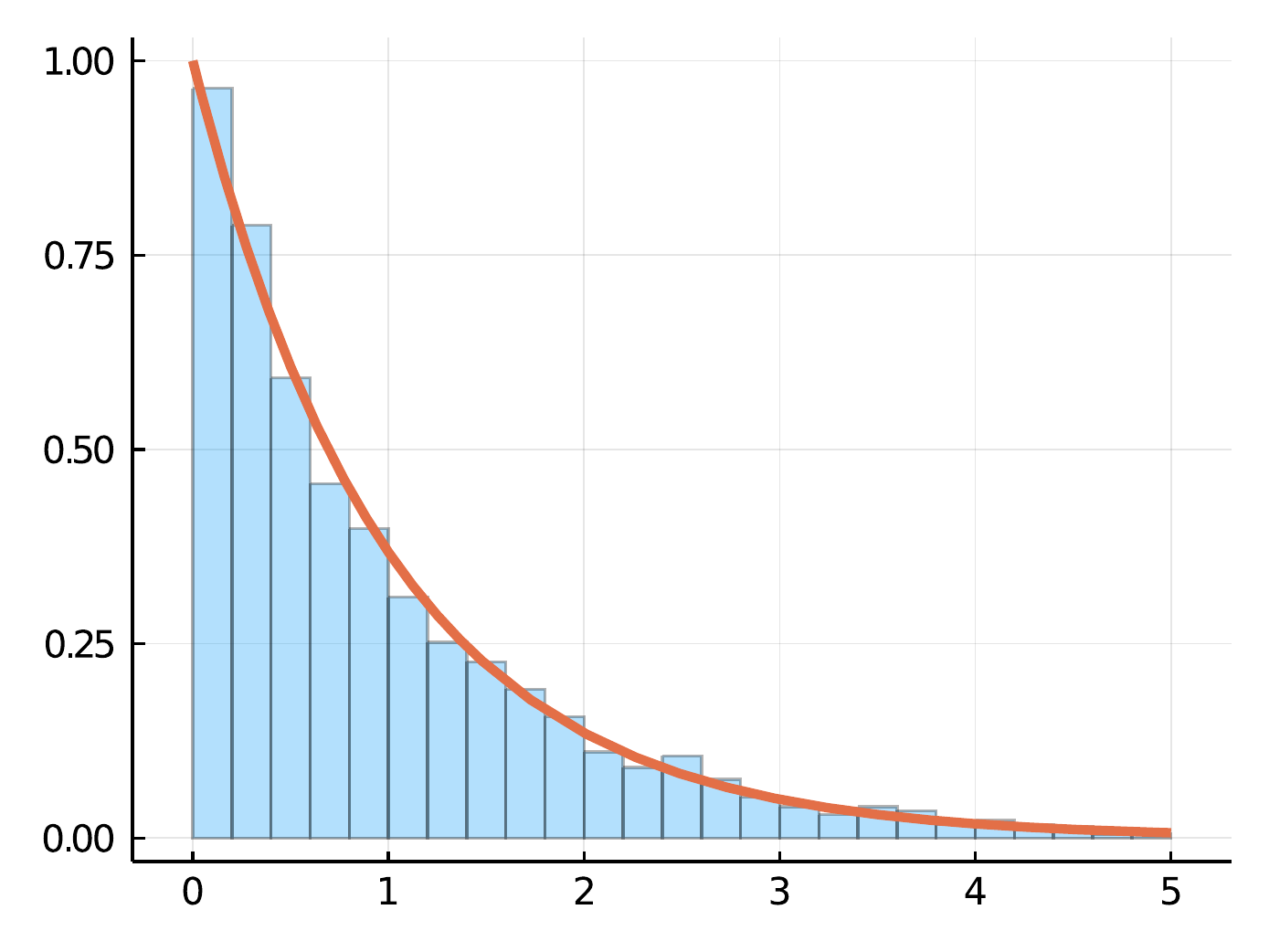}};
    \node at (1.8,1.8) {(a)};
    \node at (2.2,-1.2) {$s$};
    \node at (-1.9,1.9) {$p_0(s)$};
    
    \begin{scope}[shift = {(5.5,0)}]
    \node at (0,0) {\includegraphics[height=3.5cm]{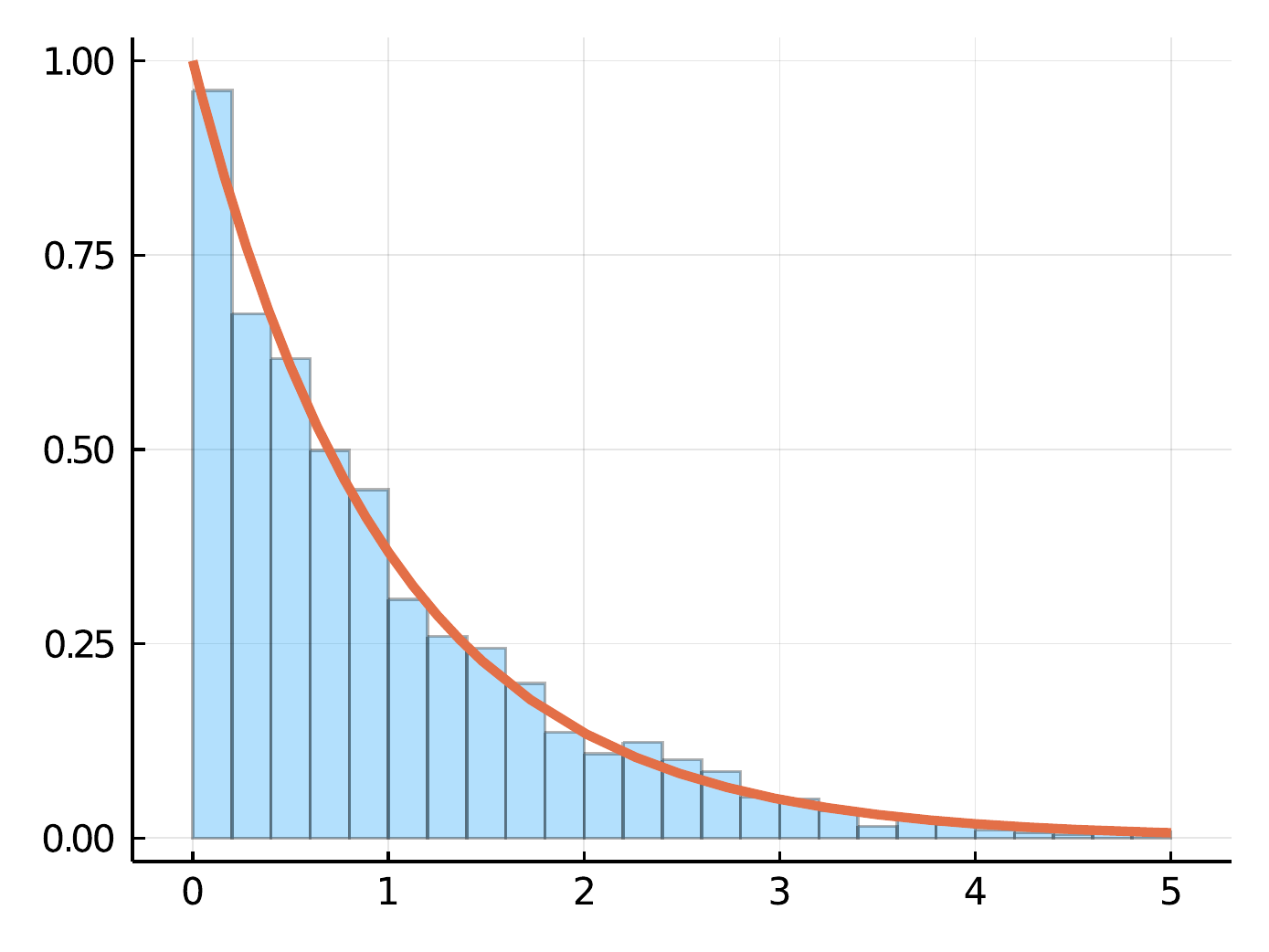}};
    \node at (1.8,1.8) {(b)};
    \node at (2.2,-1.2) {$s$};
    \node at (-1.9,1.9) {$p_0(s)$};
    \end{scope}
    
    \begin{scope}[shift = {(11,0)}]
    \node at (0,0) {\includegraphics[height=3.5cm]{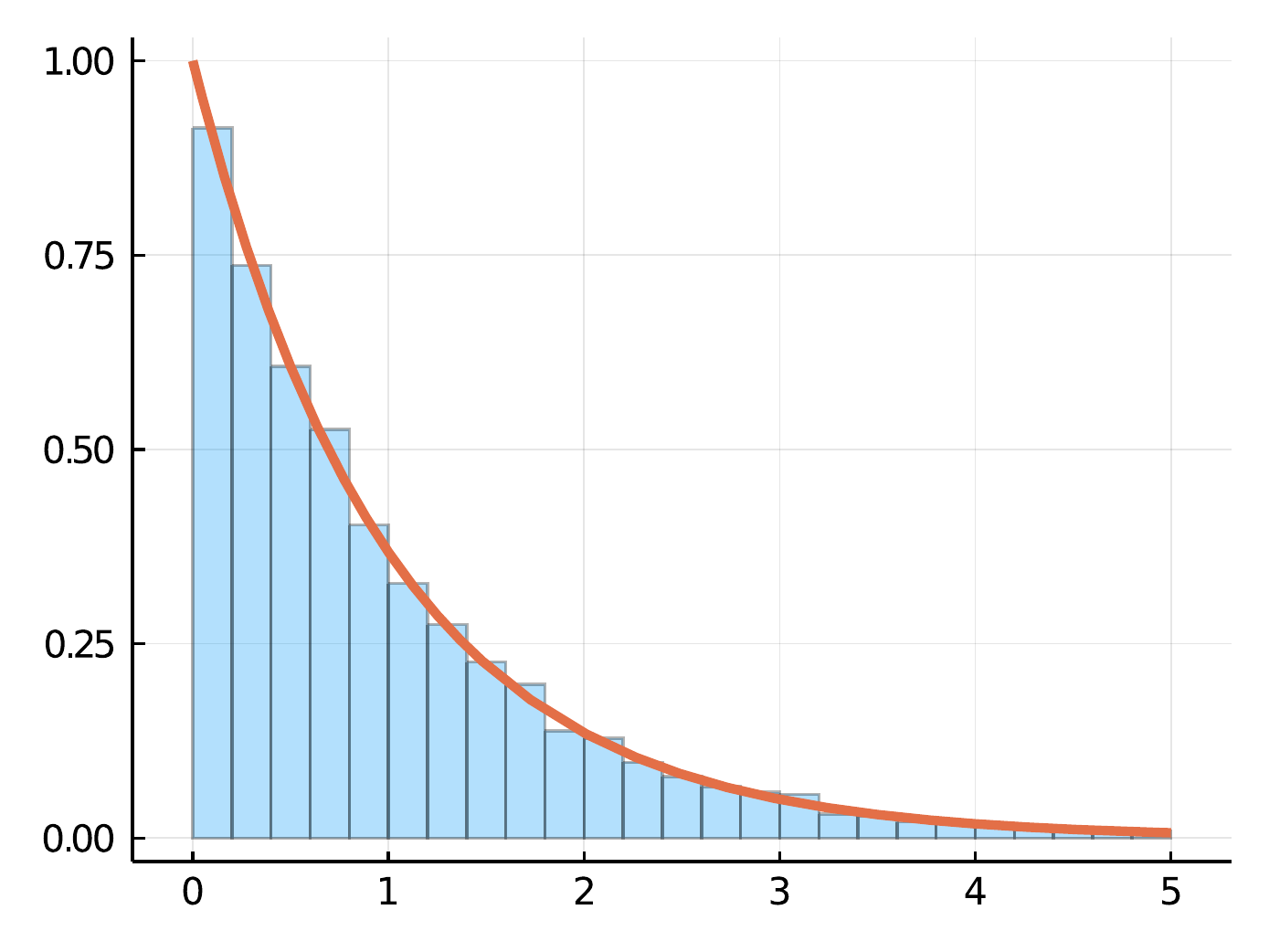}};
    \node at (1.8,1.8) {(c)};
    \node at (2.2,-1.2) {$s$};
    \node at (-1.9,1.9) {$p_0(s)$};
    \end{scope}
    \end{tikzpicture}
        
    \end{center}
    \caption{Nearest neighbor spacing distribution $p_0(s)$.
    The red line shows the curve $e^{-s}$.}
    \label{fig:sp0}
\end{figure}
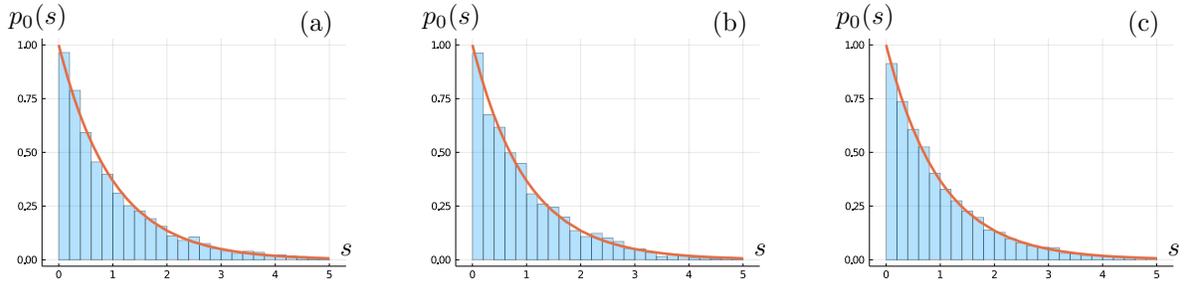

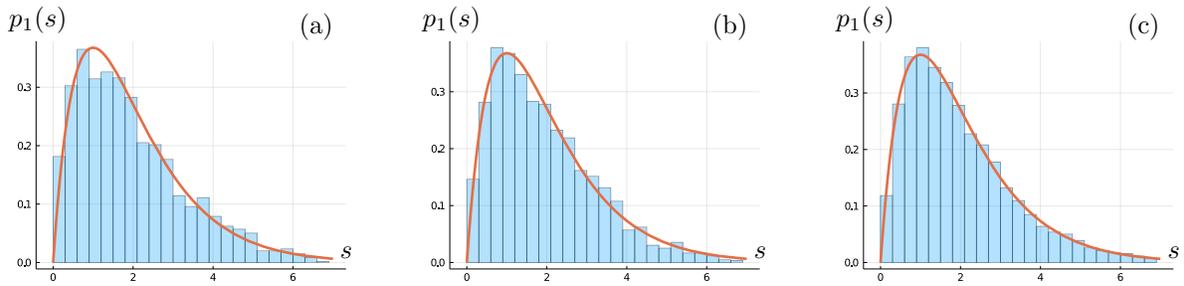
\begin{figure}
    \begin{center}
        
    \begin{tikzpicture}
    \node at (0,0) {\includegraphics[height=3.5cm]{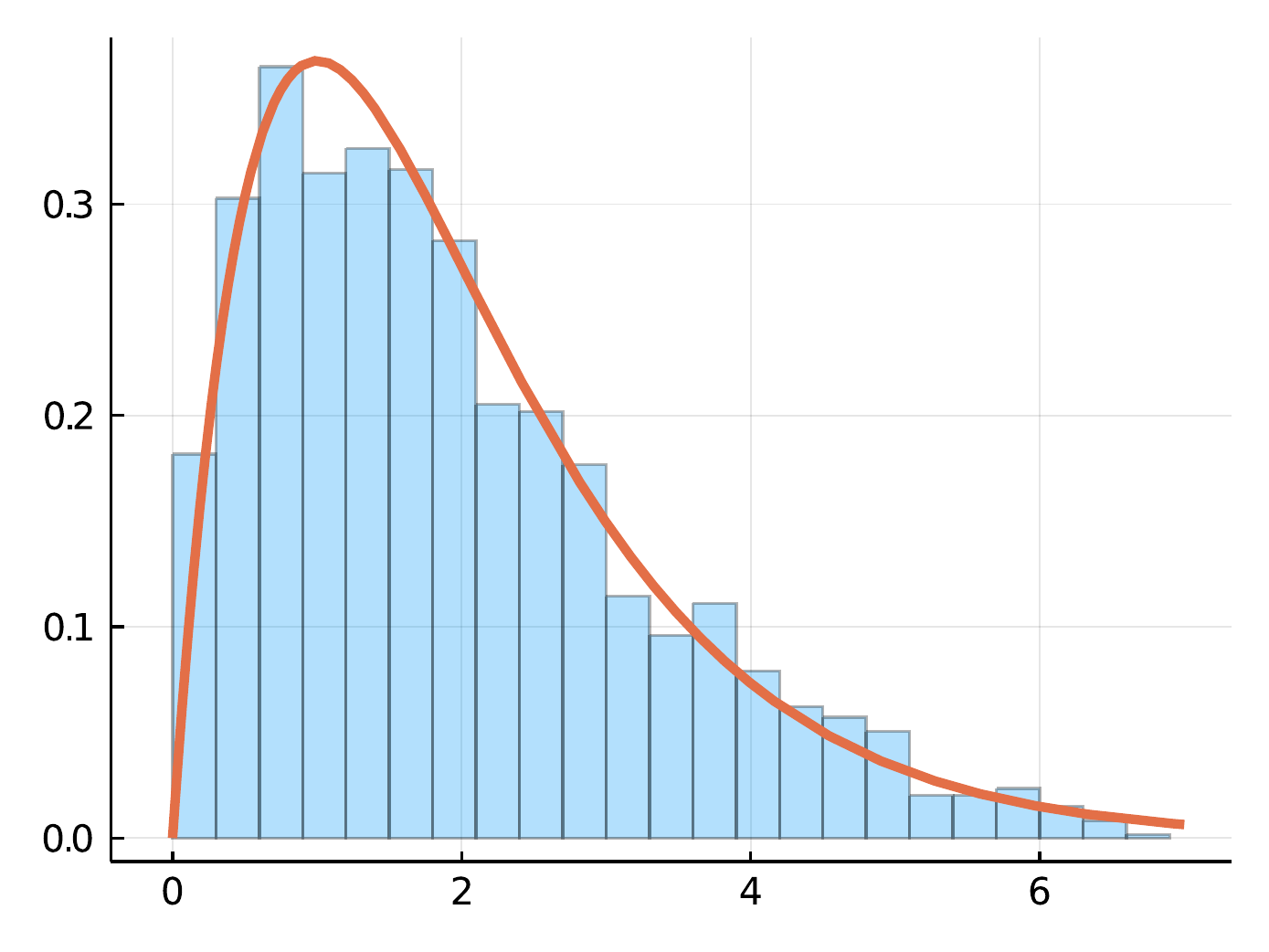}};
    \node at (1.8,1.8) {(a)};
    \node at (2.2,-1.2) {$s$};
    \node at (-1.9,1.9) {$p_1(s)$};
    
    \begin{scope}[shift = {(5.5,0)}]
    \node at (0,0) {\includegraphics[height=3.5cm]{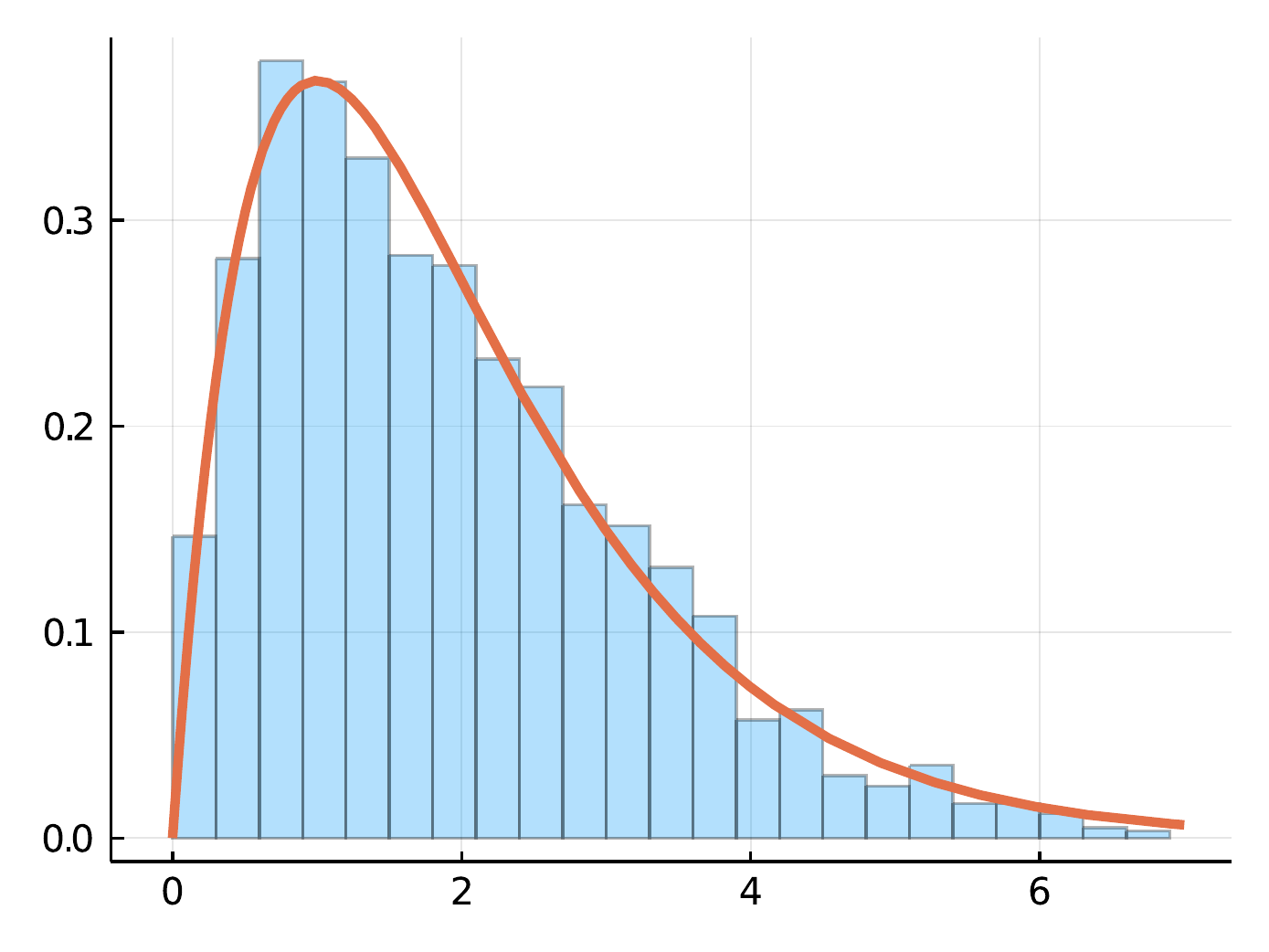}};
    \node at (1.8,1.8) {(b)};
    \node at (2.2,-1.2) {$s$};
    \node at (-1.9,1.9) {$p_1(s)$};
    \end{scope}
    
    \begin{scope}[shift = {(11,0)}]
    \node at (0,0) {\includegraphics[height=3.5cm]{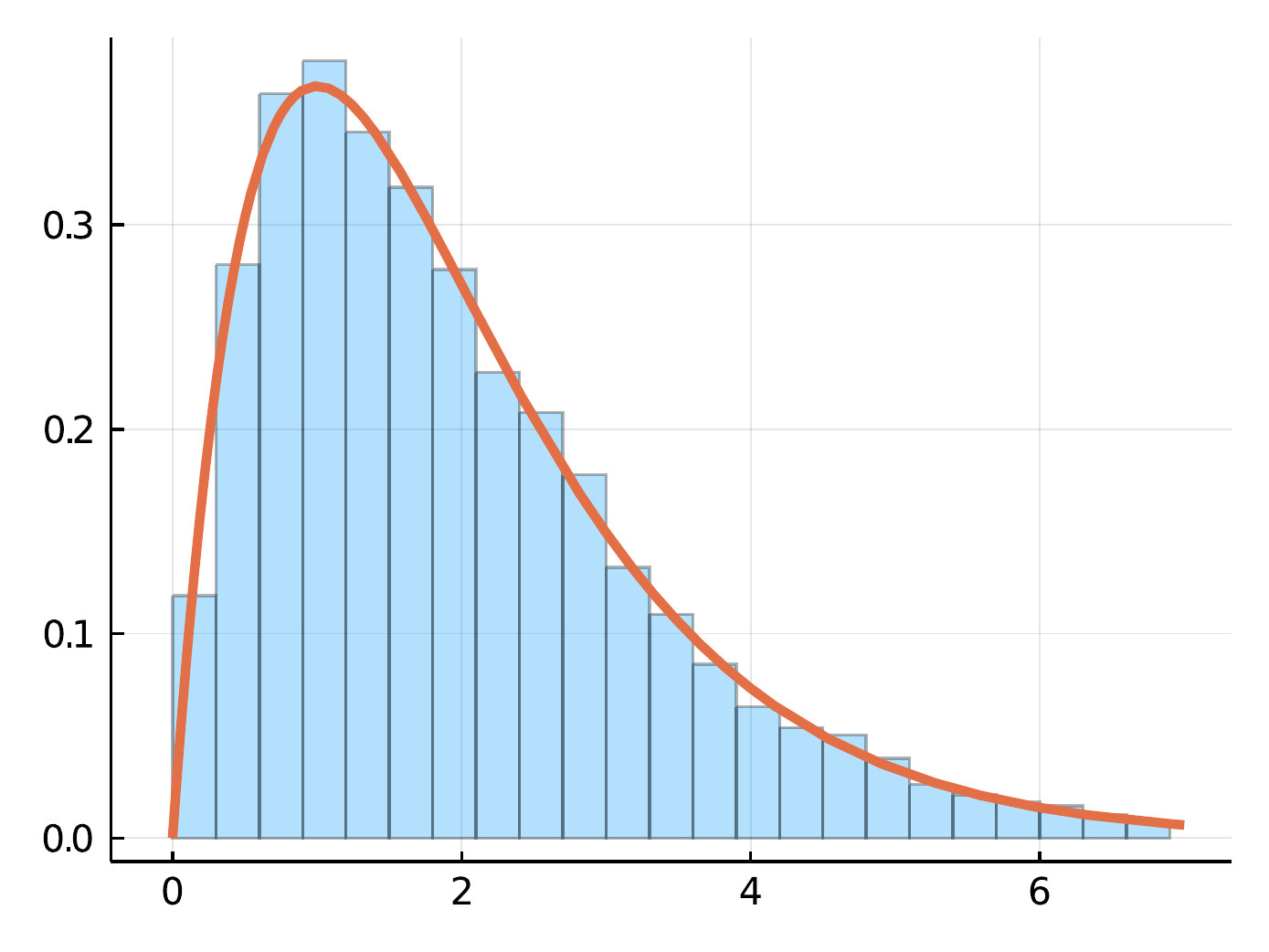}};
    \node at (1.8,1.8) {(c)};
    \node at (2.2,-1.2) {$s$};
    \node at (-1.9,1.9) {$p_1(s)$};
    \end{scope}
    \end{tikzpicture}
        
    \end{center}
    \caption{Next nearest neighbor spacing distribution $p_1(s)$.
    The red line shows the curve $s\,e^{-s}$.}
    \label{fig:sp1}
\end{figure}

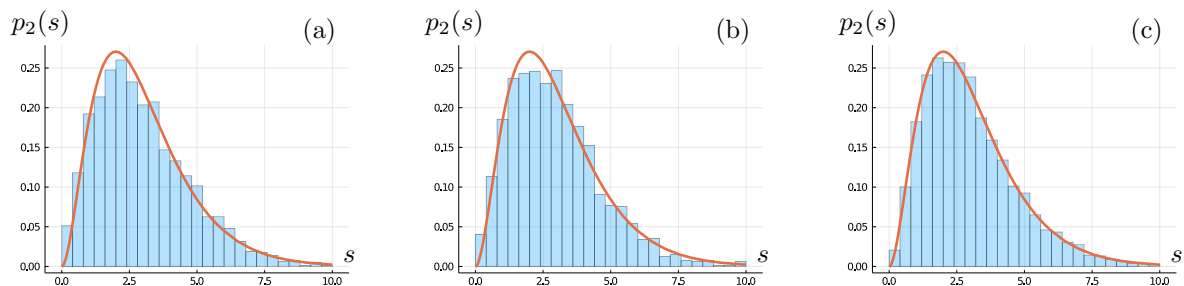
\begin{figure}
    \begin{center}
        
    \begin{tikzpicture}
    \node at (0,0) {\includegraphics[height=3.5cm]{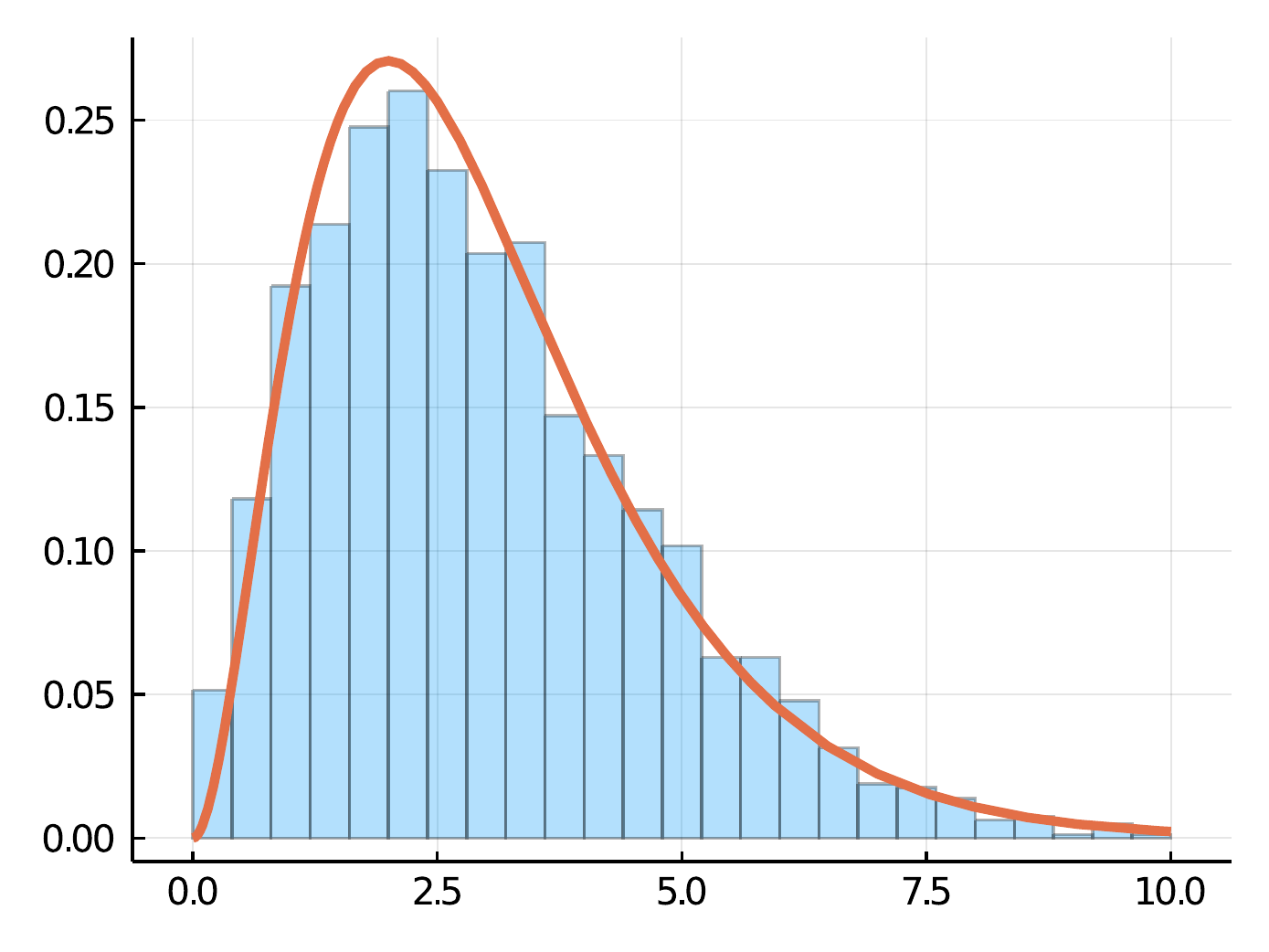}};
    \node at (1.8,1.8) {(a)};
    \node at (2.2,-1.2) {$s$};
    \node at (-1.9,1.9) {$p_2(s)$};
    
    \begin{scope}[shift = {(5.5,0)}]
    \node at (0,0) {\includegraphics[height=3.5cm]{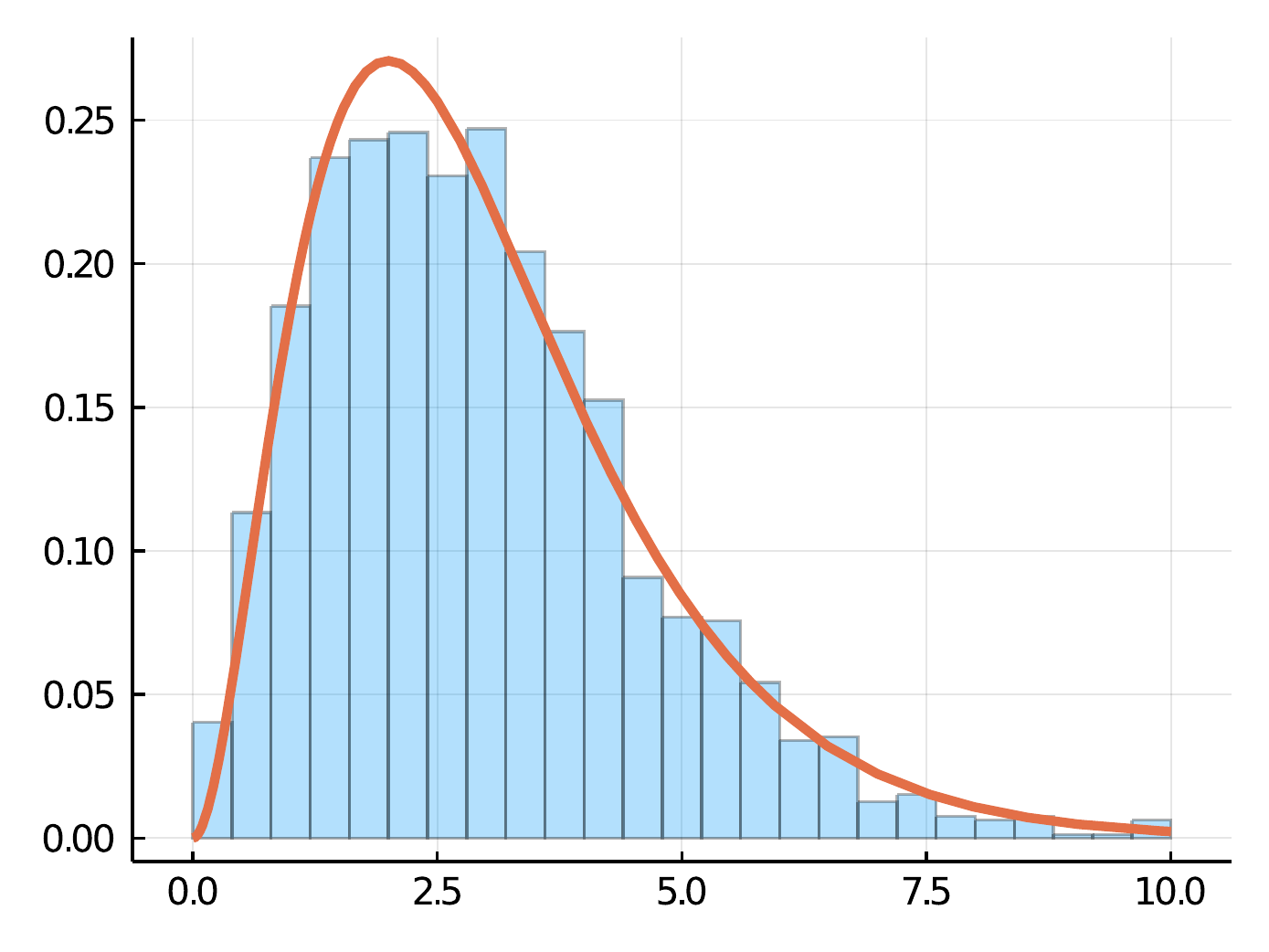}};
    \node at (1.8,1.8) {(b)};
    \node at (2.2,-1.2) {$s$};
    \node at (-1.9,1.9) {$p_2(s)$};
    \end{scope}
    
    \begin{scope}[shift = {(11,0)}]
    \node at (0,0) {\includegraphics[height=3.5cm]{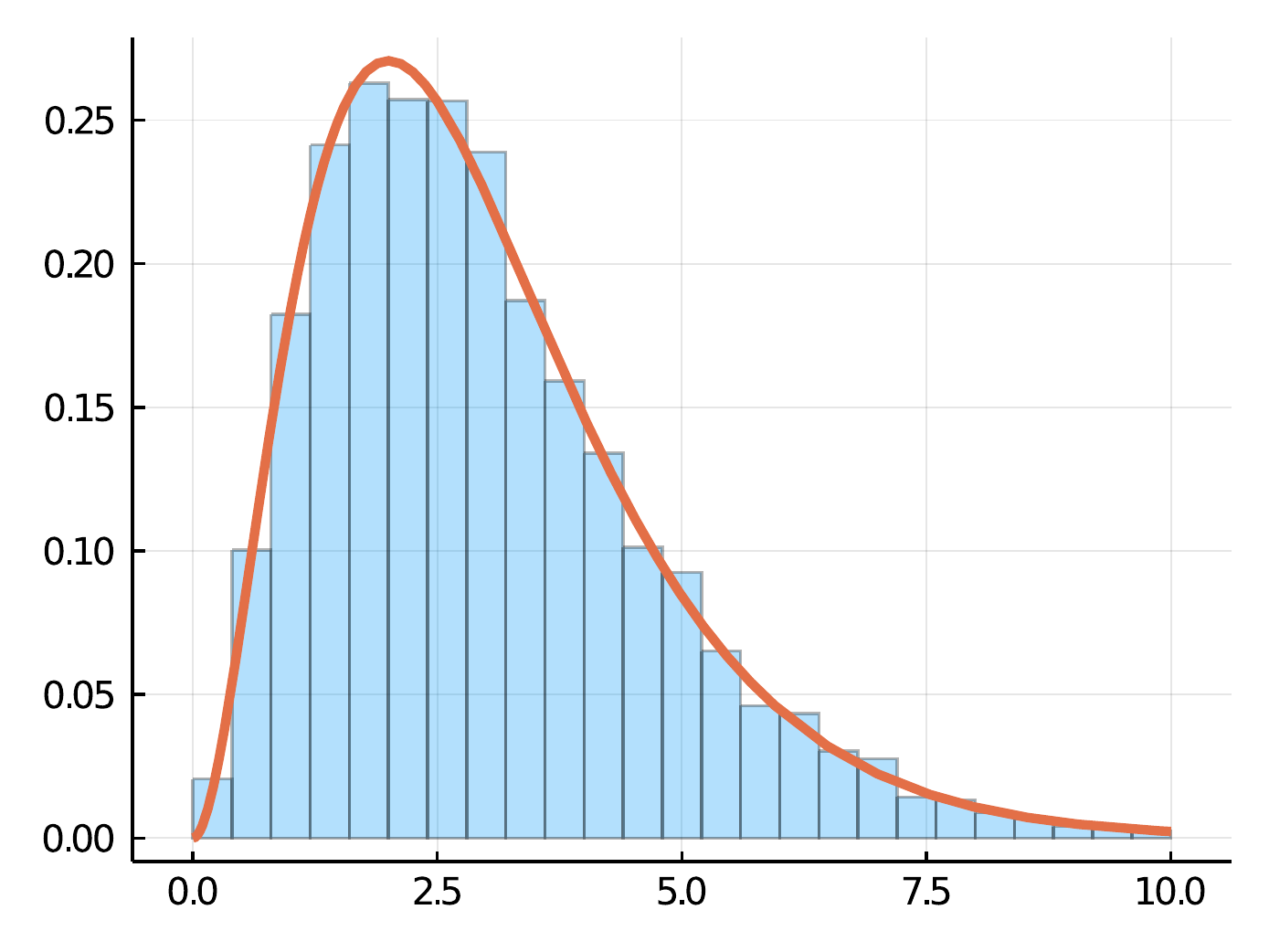}};
    \node at (1.8,1.8) {(c)};
    \node at (2.2,-1.2) {$s$};
    \node at (-1.9,1.9) {$p_2(s)$};
    \end{scope}
    \end{tikzpicture}
        
    \end{center}
    \caption{Next-next nearest neighbor spacing distribution $p_2(s)$.
    The red line shows the curve $s^2 e^{-s}/2$.}
    \label{fig:sp2}
\end{figure}

Fig.~\ref{fig:sp0}, Fig.~\ref{fig:sp1}, and Fig.~\ref{fig:sp2} show the next$^k$ nearest neighbor spacing distribution $p_k(s)$ for $k = 0, 1, 2$.
We see the agreement with the Poisson distribution $s^k e^{-s}/k!$, which provides a numerical evidence of the conjecture.

\providecommand{\bysame}{\leavevmode\hbox to3em{\hrulefill}\thinspace}
\providecommand{\MR}{\relax\ifhmode\unskip\space\fi MR }
\providecommand{\MRhref}[2]{%
  \href{http://www.ams.org/mathscinet-getitem?mr=#1}{#2}
}
\providecommand{\href}[2]{#2}


\begin{thebibliography}{HSBT10}

\bibitem[BK13]{Bourgade:2013}
P.~Bourgade and J.~P. Keating,
  \href{https://doi.org/10.1007/978-3-0348-0697-8_4}{\emph{{Quantum Chaos,
  Random Matrix Theory, and the Riemann $\zeta$-function}}}, {Chaos}, Progress
  in Mathematical Physics, vol.~66, Springer Basel, 2013, pp.~125--168.

\bibitem[CHT08]{Clozel:2008PMIHES}
L.~Clozel, M.~Harris, and R.~Taylor, \emph{{Automorphy for some $l$-adic lifts
  of automorphic mod $l$ Galois representations}},
  \href{https://doi.org/10.1007/s10240-008-0016-1}{Publ. Math. IHES
  \textbf{108} (2008)}, no.~1, 1--181.

\bibitem[HSBT10]{Harris:2010AM}
M.~Harris, N.~Shepherd-Barron, and R.~Taylor, \emph{{A family of Calabi--Yau
  varieties and potential automorphy}},
  \href{https://doi.org/10.4007/annals.2010.171.779}{Ann. Math. \textbf{171}
  (2010)}, no.~2, 779--813.

\bibitem[Kob93]{Koblitz:1993}
N.~Koblitz,
  \href{https://doi.org/10.1007/978-1-4612-0909-6}{\emph{{Introduction to
  Elliptic Curves and Modular Forms}}}, Graduate Texts in Mathematics, vol.~97,
  Springer New York, 1993.

\bibitem[KS99]{Katz:1999}
N.~M. Katz and P.~Sarnak, \emph{{Random Matrices, Frobenius Eigenvalues, and
  Monodromy}}, Colloquium Publications, vol.~45, American Mathematical Society,
  Providence, R.I, 1999.

\bibitem[Lib98]{Liboff:1998IJTP}
R.~L. Liboff, \emph{{Quasi-Chaotic Property of the Prime-Number Sequence}},
  \href{https://doi.org/10.1023/a:1026656418104}{Int. J. Theor. Phys.
  \textbf{37} (1998)}, no.~12, 3109--3117.

\bibitem[Meh04]{Mehta:2004RMT}
M.~L. Mehta, \href{https://doi.org/10.1016/S0079-8169(04)80088-6}{\emph{{Random
  Matrices}}}, 3rd ed., Pure and Applied Mathematics, vol. 142, Academic Press,
  2004.

\bibitem[Mon73]{Montgomery:1973}
H.~L. Montgomery, \emph{{The pair correlation of zeros of the zeta function}},
  Analytic number theory, Proc. Sympos. Pure Math. \textbf{XXIV} (1973),
  181--193.

\bibitem[Odl87]{Odlyzko:1987MC}
A.~M. Odlyzko, \emph{{On the distribution of spacings between zeros of the zeta
  function}}, \href{https://doi.org/10.1090/S0025-5718-1987-0866115-0}{Math.
  Comp. \textbf{48} (1987)}, 273--308.

\bibitem[Tat65]{Tate:1965}
J.~Tate, \emph{{Algebraic Cycles and Poles of Zeta Functions}}, {Arithmetical
  Algebraic Geometry}, Harper \& Row, 1965, pp.~93--110.

\bibitem[Tay08]{Taylor:2008PMIHES}
R.~Taylor, \emph{{Automorphy for some $l$-adic lifts of automorphic mod $l$
  Galois representations. {II}}},
  \href{https://doi.org/10.1007/s10240-008-0015-2}{Publ. Math. IHES
  \textbf{108} (2008)}, no.~1, 183--239.

\bibitem[Tim06]{Timberlake:2006AJP}
T.~Timberlake, \emph{{Random numbers and random matrices: Quantum chaos meets
  number theory}}, \href{https://doi.org/10.1119/1.2198883}{Amer. J. Phys.
  \textbf{74} (2006)}, no.~6, 547--553.

\bibitem[TT07]{Timberlake:2007}
T.~Timberlake and J.~Tucker, \emph{{Is there quantum chaos in the prime
  numbers?}}, \href{https://arxiv.org/abs/0708.2567}{{\ttfamily arXiv:0708.2567
  [quant-ph]}}.

\bibitem[Wol14]{Wolf:2014PRE}
M.~Wolf, \emph{{Nearest-neighbor-spacing distribution of prime numbers and
  quantum chaos}}, \href{https://doi.org/10.1103/physreve.89.022922}{Phys. Rev.
  E \textbf{89} (2014)}, no.~2, 022922,
  \href{https://arxiv.org/abs/1212.3841}{{\ttfamily arXiv:1212.3841
  [math.NT]}}.

\end{thebibliography}
\end{document}